\documentclass[11pt]{article}
\usepackage{a4,pstricks,theorem,amsmath,amssymb,latexsym,fullpage,bbm}
\usepackage{epsf,psfrag,epsfig}
\usepackage[english]{babel}
\usepackage[latin1]{inputenc}

\newcommand{\Ref}[1]{(\ref{#1})}
\newcommand{\dem}{\noindent \textbf{Proof: }} 
\newcommand{\findem}{\hfill $\square~~$\\}

\newtheorem{thm}{Theorem}

\newtheorem{lemma}[thm]{Lemma}

\newcommand{\pare}[1]{\!\!\left( #1 \right)}

\newcommand{\ite}{\noindent $\bullet$ }
\newcommand{\titre}[1]{\noindent \textbf{#1}}

\newcommand{\la}{\lambda}
\newcommand{\Om}{\Omega}
\newcommand{\eps}{\epsilon}
\newcommand{\si}{\sigma}
\newcommand{\Si}{\mathfrak{S}}
\newcommand{\mC}{\mathcal{C}}

\newcommand{\hh}{\textbf{h}}
\newcommand{\hb}{\overline{h}}
\newcommand{\hhb}{\overline{\hh}}
\newcommand{\ZZ}{\mathbb{Z}}

\newcommand{\RR}{\mathbb{R}}

\newcommand{\id}{\textrm{Id}}
\newcommand{\gr}{\textrm{gr}}
\newcommand{\Int}{\int\!\!\!\!\int}
\newcommand{\tth}{\textrm{th}}
\newcommand{\connect}{\textrm{ connected}}
\newcommand{\sign}{\textrm{sign}}
\newcommand{\Vol}{\textrm{Vol}}

\title{Solution to a combinatorial puzzle arising from Mayer's theory of cluster integrals}

\author{Olivier Bernardi}
\date{\today}

\begin{document}

\maketitle

\abstract{Mayer's theory of cluster integrals allows one to write the partition function of a gas model as a generating function of weighted graphs. Recently, Labelle, Leroux and Ducharme have studied the graph weights arising from the one-dimensional hard-core gas model and noticed that the sum of the weights over all connected graphs with $n$ vertices is $(-n)^{n-1}$. This is, up to sign, the number of rooted Cayley trees on $n$ vertices and the authors asked for a combinatorial explanation. The main goal of this article is to provide such an explanation.}

\section{Introduction}
In \cite{Mayer:theory}, Mayer used an algebraic identity in order to express the partition function of a gas model as a generating function of weighted graphs. By Mayer's transformation, any choice of an interaction potential between particles in the gas leads to a specific graph weight. For instance, in the case of the one-dimensional hard-core gas, Labelle, Leroux and Ducharme \cite{Leroux:mayers-theory2} have shown that the Mayer's weight of a connected graph $G$ having vertex set $V(G)=\{0,\ldots, n\}$ and edge set $E(G)$ is $w(G)=(-1)^{|E(G)|}\Vol(\Pi_G)$ where $\Vol(\Pi_G)$ is the volume of the {$n$-dimensional} polytope 
$$\Pi_G=\{(x_1,\ldots,x_{n})\in \RR^n/~x_0=0 \textrm{ and } |x_i-x_j|\leq 1 \textrm{ for all edge } (i,j)\in E(G)\}.$$
The pressure in the model is related to Mayer's weights by 
\begin{eqnarray}\label{eq:pressure-Mayer-intro}
P~=~kT \sum_{G \textrm{ connected graph}}  w(G)\frac{z^{|V(G)|}}{|V(G)|!},
\end{eqnarray}
where $k$ is Boltzmann's constant, $T$ is the temperature and $z$ is the \emph{activity}. \\

It is known (see \cite{Bridge:Dimensional-reduction}) that the pressure of the hard-core gas is $P=kT L(z)$, where $L(z)$ is the \emph{Lambert function} defined by the functional equation $L(z)=z\exp(-L(z))$. 
Comparing this expression of the pressure with \Ref{eq:pressure-Mayer-intro} and extracting the coefficient of $z^{n+1}$ gives
\begin{eqnarray} \label{eq:Mayer-continuum-gas-intro}
\sum_{G\in \mC_n}w(G)~=~(-1)^n(n+1)^{n},
\end{eqnarray}
where the sum is over all connected graphs with $n+1$ vertices.   Labelle \emph{et al.} observed that the right-hand-side of~\Ref{eq:Mayer-continuum-gas-intro} is, up to sign, the number of rooted Cayley trees with $n+1$ vertices and asked for a combinatorial explanation \cite[Question 1]{Leroux:mayers-theory2}. The main purpose of this paper is to give such an explanation.\\

The outline of the paper is as follows. In Section~\ref{section:review-Mayer}, we briefly review  Mayer's theory of cluster integrals following the line of \cite{Leroux:mayers-theory1}. We illustrate this theory on a very simple model of \emph{discrete gas} and prove the equivalence of this model with the Potts model on the complete graph. Comparing two expressions of the pressure in the discrete gas leads to a surprising combinatorial identity.  In Section~\ref{section:discrete-gas}, we give a combinatorial proof of this identity. 
In Section~\ref{section:continuum-gas}, we recall Mayer's setting for the hard-core continuum gas and then give a combinatorial proof of Equation~\Ref{eq:Mayer-continuum-gas-intro}, thereby answering the question of Labelle \emph{et al}.\\

We close this section with some notations. We denote by  $\mathbb{Z}$ the set of integers and by $\RR$ the set of real numbers. We denote $[n]=\{1,\ldots,n\}$ and by $\Si_n$ the set of permutations of $[n]$. 
In this paper, all \emph{graphs} are \emph{simple}, \emph{undirected} and \emph{labelled}. Let $G$ be a graph. 
 We denote by $v(G)$, $e(G)$ and $c(G)$ respectively the number of vertices, edges and connected components of $G$. A graph $H$ is a \emph{spanning subgraph} of $G$ if the vertex sets of $H$ and $G$ are the same while the edge set of $H$ is included in the edge set of $G$; we denote $H\subseteq G$ in this case. We denote by $e=(i,j)$ the edge with endpoints $i$ an $j$ and write $e\in G$ if the edge $e$ belongs to $G$. For any edge $e$, we denote by $G\oplus e$ the graph obtained from $G$ by either adding the edge $e$ if $e\notin G$ or by deleting this edge if $e\in G$. \\

\section{Review of Mayer's theory of cluster integrals}\label{section:review-Mayer}
Consider a gas made of $n$ (indistinguishable) particles in a vessel $\Om\subset\RR^d$. We suppose that the gas is free from outside influence and that interaction between two particles $i$ and $j$ at positions $x_i$ and $x_j$ is given by the potential $\phi(x_i,x_j)$. In the classical Boltzmann setting, the probability measure of a configuration is proportional to $\exp(-H/kT)$, where  $k$ is Boltzmann's constant, $T$ is the temperature and $H$ is the \emph{Hamiltonian} of the system given by
$$H=\sum_{1\leq i\leq n}\frac{m_i v_i^2}{2}+\sum_{1\leq i<j\leq n}\phi(x_i,x_j),$$
where $x_i,v_i,m_i$  and  $\frac{m_i v_i^2}{2}$ are respectively the position, velocity, mass and kinetic energy of the $i^\tth$ particle.\\

The \emph{partition function} of the gas model is 
$$
Z(\Omega,T,n) =\frac{1}{h^{dn}n!} \Int_{x_1,\ldots,x_n \in\Omega,~v_1,\ldots,v_n\in\RR^d}exp(-H/kT)dx_1\ldots dx_ndv_1\ldots dv_n, 
$$
where $h$ is Planck's constant. After integrating over all possible velocity, the partition function becomes
$$
Z(\Omega,T,n) =\frac{1}{\lambda^n n!} \Int_{\Omega^n} \prod_{i<j} \exp\pare{-\frac{\phi(x_i,x_j)}{kT}}dx_1\ldots dx_n,
$$
where  $\lambda$ depends on the temperature $T$. \\

Mayer noticed that the partition function can be decomposed into a sum over graphs. Indeed, by setting 
$f(x_i,x_j)=\exp\pare{-\frac{\phi(x_i,x_j)}{kT}}-1,$ one gets
$$\prod_{i<j} \exp\pare{-\frac{\phi(x_i,x_j)}{kT}}~=~\prod_{i<j} 1+f(x_i,x_j)~=~\sum_{G\subseteq K_n} \prod_{(i,j)\in G}f(x_i,x_j),$$
where the sum is over all graphs on $n$ vertices (equivalently, spanning subgraph of the complete graph $K_n$) and the inner product is over all edges of $G$. In terms of the partition function, this gives \emph{Mayer's relation}:
$$
\lambda^n n!Z(\Omega,T,n)~=~ \Int_{\Omega^n} \prod_{i<j} \exp\pare{-\frac{\phi(x_i,x_j)}{kT}}dx_1\ldots dx_n~=~ \sum_{G\subseteq K_n} W(G), \\
$$
where  $\displaystyle ~W(G)=\Int_{\Omega^n}\prod_{(i,j)\in G}f(x_i,x_j)dx_1\ldots dx_n~$ is the \emph{first Mayer's weight} of the graph $G$.\\

\bigskip

\titre{Example: the \emph{discrete gas}.} Suppose $\Omega$ is made of $q$ distinct boxes $B_1,\ldots,B_q$ of volume 1 and that the interaction potential $\phi(x_i,x_j)$ is equal to $\alpha$ if the particles $i$ and $j$ are in the same box and 0 otherwise. 
By definition, the Mayer's weight of a graph $G$ is 
$$W(G)=\Int_{\Omega^n}\prod_{(i,j)\in G}f(x_i,x_j)dx_1\ldots dx_n=\sum_{c:[n]\mapsto [q]} \Int_{x_1\in B_{c(1)},\ldots, x_n\in B_{c(n)}}\prod_{i<j} f(x_i,x_j)dx_1\ldots dx_n.$$
We denote $u=\exp(-\alpha/kT)$ and observe that $f(x_i,x_j)=u-1$ if $i$ and $j$ are in the same box and 0 otherwise. 
Therefore, the product $\prod_{(i,j)\in G}f(x_i,x_j)$ equals  $(u-1)^{e(G)}$ if the value of $c$ is constant over each connected components of the graph $G$ and 0 otherwise. Summing over all mappings $c:[n]\mapsto [q]$ gives 
$$W(G)~=~ q^{c(G)} (u-1)^{e(G)},$$
since there are~$q^{c(G)}$ mappings  $c:[n]\mapsto [q]$ which are constant over each connected components of~$G$.\\

In our discrete gas example, a direct calculation of the partition function  gives 
$$
\lambda^n n! Z(\Om,T,n) = \sum_{c:[n]\mapsto [q]} \Int_{x_1\in B_{c(1)}\ldots x_n\in B_{c(n)}}\prod_{i<j} \exp\pare{-\frac{\phi(x_i,x_j)}{kT}}dx_1\ldots dx_n = \sum_{c:[n]\mapsto [q]}u^{\delta(c)}, 
$$
where $\delta(c)$ is the number of edges $(i,j)\in K_n$ such that $c(i)=c(j)$. Hence, Mayer's relation reads 
\begin{eqnarray} \label{eq:Tutte=Potts}
\sum_{c:[n]\mapsto [q]}u^{\delta(c)}=\sum_{G\subseteq K_n} q^{c(G)}(u-1)^{e(G)}.
\end{eqnarray}
Equation \Ref{eq:Tutte=Potts} is a special case of the equivalence established by Fortuin and Kastelein \cite{Fortuin:Tutte=Potts} between the partition function of the Potts model  (see~e.g.~\cite{Baxter:exactly-solved-model}) and the \emph{Tutte polynomial} (see~e.g.~\cite{Bollobas:Tutte-poly}). Indeed, the right-hand-side corresponds to the partition function of the Potts model on the complete graph $K_n$ while the left-hand-side corresponds to the subgraph expansion of the Tutte polynomial of $K_n$ up to scaling and change of variables.  The relation of Fortuin and Kastelein is the generalisation of~\Ref{eq:Tutte=Potts} obtained by replacing the complete graph~$K_n$ by any graph~$H$. This more general case relies on the observation that $\prod_{(i,j)\in H} \phi_{i,j} = \sum_{G\subseteq H} \prod_{(i,j)\in G} f_{i,j} $ 
as soon as $\phi_{i,j}=1+f_{i,j} $ for all $(i,j)\in H$.\\

\bigskip

We now return to the general theory of Mayer and consider a system with an arbitrary number of particles. The \emph{grand canonical partition function} is defined by 
$$Z_\gr(z)~\equiv~ Z_\gr(\Om,T,z)~=~\sum_{n\leq 0} z^n \la^n Z(\Om,T,n),$$
where $z$ is the \emph{activity} of the system.  In terms of Mayer's weights, the grand canonical partition function is the exponential generating functions of graphs weighted by their first Mayer's weight:
$$Z_\gr(z)~=~\sum_{n\leq 0} \frac{z^n}{n!} \sum_{G\subseteq K_n} W(G)~=~ \sum_{G} W(G)\frac{z^{v(G)}}{v(G)!}.$$
 
The macroscopic parameters of the systems, such as the density $\rho$, or pressure $P$, can be obtained from $Z_\gr(z)$ by the relations 
$$P=\frac{kT}{|\Om|}\log(Z_\gr(z)) ~~\textrm{ and }~~ \rho=\frac{z}{|\Om|}\frac{\partial}{\partial z} \log(Z_\gr(z)).$$

Observe that the first Mayer's weight is \emph{multiplicative} over connected components, that is, if a graph $G$ is the disjoint union of two graphs $G_1$ and $G_2$ then
$W(G)=W(G_1)W(G_2)$. This is the key property implying  
$ ~\log(Z_\gr(z))=\sum_{G \connect} W(G)z^{v(G)}~$ (see \cite{Leroux:mayers-theory2} for a complete proof), or equivalently,  
\begin{eqnarray}\label{eq:pressure}
P=\frac{kT}{|\Om|}\sum_{G \connect} W(G)\frac{z^{v(G)}}{v(G)!}~.
\end{eqnarray}

\bigskip

\titre{Example: the discrete gas.} 
For the discrete gas model introduced before, Equation~\Ref{eq:pressure} gives 
\begin{eqnarray}\label{eq:pressure-discrete}
\frac{P}{kT}= \sum_{G \connect} (u-1)^{e(G)} \frac{z^{v(G)}}{v(G)!}.
\end{eqnarray}
In the special case of an infinite repulsive interaction between particles in the same box,  that is, $\alpha=\infty$ and $u=0$, the pressure $P$ can also be computed directly. Indeed, in this case, one gets
$$\lambda^n n! Z(\Om,T,n)~=~\sum_{c:[n]\mapsto [q]}u^{\delta(c)}~=~\#\{c:[n]\mapsto [q] ~\textrm{injective}\} ~=~q(q-1)\cdots(q-n+1),$$
and 
$$Z_\gr(z)~=~\sum_{n\geq 0} z^n \la^n Z(\Om,T,n)~=~\sum_{n\geq 0}{q \choose n}z^n~=~(1+z)^q.$$
This expression for the grand canonical partition function comes to no surprise since each of the $q$ boxes contains either nothing (activity 1) or one particle (activity $z$). Now, 
$$\frac{P}{kT}~=~\frac{1}{q}\log(Z_\gr(z))~=~\log(1+z)~=~\sum_{n>0}\frac{(-1)^{n-1}}{n}z^n,$$
and extracting the coefficient of $z^n$ in both side of \Ref{eq:pressure-discrete} gives
\begin{eqnarray}\label{eq:pressure-discrete-coeff}
(-1)^{n-1} \sum_{G\subseteq K_n \connect} (-1)^{e(G)}~=~(n-1)!.
\end{eqnarray}

Identity~\Ref{eq:pressure-discrete-coeff} is quite surprising at first sight but can be understood by recognising in the left-hand-side the evaluation of the Tutte polynomial of $K_n$ counting the root-connected acyclic orientations (acyclic orientation in which the vertex 1 is the only source) \cite{Greene:interpretation-Tutte-poly}. Indeed, in the case of the complete graph $K_n$, root-connected acyclic orientations are linear orderings of $[n]$ in which $1$ is the least element, or equivalently, permutations of  $\{2,\ldots,n\}$. In the next section, we give a combinatorial proof of Equation~\Ref{eq:pressure-discrete-coeff} which avoids introducing the whole theory of the Tutte polynomial (though it is based on it) and prepares for the more evolved proof of Equation~\Ref{eq:Mayer-continuum-gas-intro}.

\section{Pressure in the hard-core discrete gas and increasing trees}\label{section:discrete-gas}
In this section, we give a combinatorial proof of \Ref{eq:pressure-discrete-coeff} by exhibiting an involution $\Phi$ on connected graphs which cancels the contribution of almost all graphs in the sum $\displaystyle \sum_{G\subseteq K_n \connect} (-1)^{e(G)}$.\\

We consider the \emph{lexicographic order} on the edges of $K_n$ defined by $(i,j)<(k,l)$ if either $\min(i,j)<\min(k,l)$ or $\min(i,j)=\min(k,l)$ and $\max(i,j)<\max(k,l)$. For a graph $G$ and an edge $e=(i,j)$ (not necessarily in $G$), we denote by  $G^{>e}$ the spanning subgraph of $G$ made of the edges which are greater than $e$.  We say that  $e=(i,j)$ is \emph{$G$-active} if there is a path in $G^{>e}$ connecting $i$ and $j$ and we denote by $e^*_{G}$ the least $G$-active edge (if there are some). We then define a mapping $\Psi$ on the set of connected graphs by setting: $\Psi(G)=G$ if there is no $G$-active edge and $\Psi(G)=G\oplus e^*_{G}$  otherwise.

\begin{lemma}\label{lem:killing-involution}
The mapping $\Psi$ is an involution on connected graphs. 
\end{lemma}

\dem 
\ite First observe that \emph{the image of a connected graph is connected}. Indeed, if the edge $e^*_{G}$ exists and belongs to $G$, then it is in a cycle of $G$ and deleting it does not disconnect $G$. \\
\ite We now prove that \emph{any edge is $G$-active if and only if it is $\Psi(G)$-active}. Suppose that the edge $e=(i,j)$ is $G$-active and let $P$ be a path of $G^{>e}$ connecting $i$ and $j$. Since $e^*_{G}\leq e$ the path $P$ does not contain $e^*_{G}$, hence $P\subseteq \Psi(G)^{>e}$ and $e$ is $\Psi(G)$-active. 
Suppose conversely that $e=(i,j)$ is  $\Psi(G)$-active and let $P$ be a path of $\Psi(G)^{>e}$ connecting $i$ and $j$. If $P$ does not contain $e^*_{G}$, then $P\subseteq  G^{>e}$ and $e$ is $G$-active. Otherwise, $e^*_{G}>e$ and there is a path  $Q$ of  $G^{>e^*_{G}}\subseteq G^{>e}$ connecting the endpoints of $e^*_{G}$. Thus, there is a path contained in $(P-e^*_{G})\cup Q\subseteq G^{>e} $ connecting $i$ and $j$ and again $e$ is $G$-active.\\
\ite By the preceding point, there is a $G$-active edge if and only if there is a $\Psi(G)$-active edge and in this case   $e^*_{\Psi(G)}=e^*_{G}$. Thus, $\Psi(\Psi(G))=G$.
\findem

The mapping $\Psi$ is an involution and $(-1)^{e(G)}+(-1)^{e(\Psi(G))}=0$ whenever $G\neq \Psi(G)$, hence
\begin{eqnarray}\label{eq:killing-discrete}
\sum_{G\subseteq K_n \connect \atop ~} (-1)^{e(G)}=\sum_{G\subseteq K_n \connect, \Psi(G)=G}(-1)^{e(G)}.
\end{eqnarray}
We now characterise the fixed points of the involution $\Psi$. 
A tree on $\{1,\ldots,n\}$ is said \emph{increasing} if the labels of the vertices are  increasing along any simple path starting from the vertex $1$.
\begin{lemma}\label{lem:survivors}
A connected graph $G$ has no $G$-active edge if and only if it is an increasing tree.
\end{lemma}

\dem
\ite We suppose first that $G$ is an increasing tree and want to prove that no edge is $G$-active. Since $G$ has no cycle, no edge in $G$ is $G$-active. Consider now an edge $e=(i,j)\notin G$ and the nearest common ancestor $k$ of $i$ and $j$ (the root vertex of $G$ being the vertex 1).  There is an edge $e'=(k,l)$ containing $k$ on the path of $G$ connecting $i$ and $j$. Since $G$ is an increasing tree, $k\leq \min(i,j)$ and $l\leq \max(i,j)$. Thus, $e'=(k,l)<e=(i,j)$ and $e$ is not $G$-active.\\
\ite Suppose now that there is no $G$-active edge. First observe that $G$ is a tree since if $G$ had a cycle then the minimal edge in this cycle would be active. We now want to prove that the tree $G$ is increasing. Suppose the contrary and consider a sequence of labels $1=i_1<i_2<\ldots<i_r>i_{r+1}$ on a path of $G$ starting from the vertex $i_1=1$. Then, the edge $(i_{r-1},i_{r+1})$ is $G$-active and we reach a contradiction.
\findem

By Lemma~\ref{lem:survivors}, the fixed points of the involution $\Psi$ are the increasing trees. The increasing trees on $\{1,\ldots,n\}$ are known to be in bijection with the permutations of $\{2,\ldots,n\}$ \cite{Stanley:volume1}. Hence, there are  $(n-1)!$ increasing trees on $[n]$ and continuing Equation \Ref{eq:killing-discrete} gives 
$$(-1)^{n-1}\!\!\!\!\sum_{G\subseteq K_n \connect}\!\!\!\!\!\!\!\! (-1)^{e(G)}=\sum_{G\subseteq K_n\connect\atop  \Psi(G)=G}\!\!\!\!\!\!\!\!\!\!\!\!(-1)^{e(G)+n-1}=\# \{ \textrm{increasing trees on } [n]\}=(n-1)!.$$
This completes the proof of Equation~\Ref{eq:pressure-discrete-coeff}.

\section{Pressure in the hard-core continuum gas and Cayley trees}\label{section:continuum-gas}
In the 1-dimensional hard-core continuum gas, the vessel is an interval $\Om=[-q/2,q/2]$ and the potential of interaction between two particles $i$ and $j$ is $\phi(x_i,x_j)=\infty$ if $|x_i-x_j|\leq 1$ and 0 otherwise. By definition, $f(x_i,x_j)\equiv \exp\pare{-\phi(x_i,x_j)/kT}-1$ is equal to -1 if $|x_i-x_j|\leq 1$ and 0 otherwise. Thus, the first Mayer weight of a graph $G$ on $n$ vertices, is
$$W(q,T,G)=\Int_{[-\frac{q}{2},\frac{q}{2}]^n} \prod_{(i,j)\in G} f(x_i,x_j) dx_1\ldots dx_n=(-1)^{e(G)} \Int_{[-\frac{q}{2},\frac{q}{2}]^n}\prod_{(i,j)\in G}\mathbbm{1}_{|x_i-x_j|\leq 1} .$$

\begin{figure}[ht!]\begin{center} \input{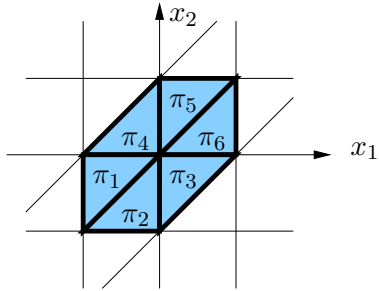}\caption{The polytope $\Pi_{K_3}$ (dashed) and its decomposition into subpolytopes.}\label{fig:polytope-sliced} \end{center}\end{figure}

In the thermodynamic limit where the volume $q$ of the Vessel tends to infinity, it becomes interesting to consider the \emph{second Mayer's weight} of connected graphs defined by $w(G)=\lim_{q\to \infty} \frac{W(q,T,G)}{q}$ and related to the pressure by $P=kT\sum_{G \connect}w(G)z^G$. In \cite{Leroux:mayers-theory2}, it is shown that for any connected graph $G$ on $\{0,\ldots,n\}$, the second Mayer weight $w(G)$ equals $(-1)^{e(G)} \Vol(\Pi_G)$, where $\Vol(\Pi_G)$ is the volume of the $n$-dimensional polytope  
$$ \Pi_G=\{(x_1,\ldots,x_n)\in \RR^n/~ x_0=0 \textrm{ and } |x_i-x_j|\le 1 \textrm{ for all edges } (i,j)\in G\}.$$
For instance, the polytope $\Pi_{K_3}$ is represented in Figure~\ref{fig:polytope-sliced}. The rest of this paper is devoted to the proof of Equation~\Ref{eq:Mayer-continuum-gas-intro} given in the introduction.

\titre{Subpolytopes.}
As observed by Bodo Lass \cite{Lass:subpolytopes-decomposition}, it is possible to decompose the polytope $\Pi_G$ into subpolytopes of volume $1/n!$. Each subpolytope is defined by fixing the integral parts and the relative order of the fractional parts of the coordinates  $x_1,\ldots,x_n$. Let us first recall some definitions: for any real number $x$ we write $x=h(x)+\eps(x)$ where $h(x)\in \ZZ$ is the \emph{integral part} and $0\le \eps(x)< 1$ is the \emph{fractional part}. Given a vector of integers $\hh=(h_1,\ldots,h_n)\in \ZZ^n$ and a permutation $\sigma\in \Si_n$, we denote by $\pi(\hh,\sigma)$ the polytope whose interior is made of the points $(x_1,\ldots,x_n)$ such that $h(x_i)=h_i$ for all $i=1\ldots n$ and  $0<\eps(x_{\si^{-1}(1)})< \cdots <\eps(x_{\si^{-1}(n)})<1$. In particular, the polytope $\pi(\hh,\sigma)$ contains the point $(h_1+\frac{\si(1)}{n+1},\ldots,h_n+\frac{\si(n)}{n+1},)$ in its interior.\\

Observe that the condition $|x_i-x_j|<1$ is equivalent to $h(x_i)-h(x_j)\in \{0,\sign(\eps(x_j)-\eps(x_i))\}$ where the value of $\sign(x)$ is -1 if $x<0$, +1 if $x>0$ and 0 if $x=0$. Therefore, a point $(x_1,\ldots,x_n)$ in the interior  of the polytope $\pi(\hh,\sigma)$ is in $\Pi_G$ if and only if  $(i,j)\in G$ implies  $h_i-h_j\in \{0, \sign(\si(j)-\si(i))\}$ with the convention that $h_0=0$ and $\si(0)=0$. This condition only depends on the pair $(\hh,\sigma)$, therefore either the polytope $\pi(\hh,\sigma)$ is included in the polytope $\Pi_G$ or the interiors of the two polytopes are disjoints. Moreover, a simple calculation shows that the volume of the polytope $\pi(\hh,\sigma)$ is $1/n!$ for all $\hh\in \ZZ^n,\sigma\in \Si_n$. This proves the following lemma. 
\begin{lemma}\label{lem:Lass}
For any connected graph $G$ on $\{0,\ldots,n\}$, the value $n!\Vol(\Pi_G)$ counts the pairs $\hh\in \ZZ^n,\si\in\Si_n$ such that  $\pi(\hh,\sigma)$ is a subpolytope of $\Pi_G$.
\end{lemma}

For example, the polytope $\Pi_{K_3}$ represented in Figure~\ref{fig:polytope-sliced} contains 6 subpolytopes $\pi_1=\pi((-1,-1),12)$,  $\pi_2=\pi((-1,-1),21)$, $\pi_3=\pi((0,-1),12)$, $\pi_4=\pi((-1,0),21)$, $\pi_5=\pi((0,0),12)$, $\pi_6=\pi((0,0),21)$ each having volume $1/2$.\\

\bigskip

\titre{Rearrangement.} 
Summing the second Mayer's weights over connected graphs and using Lemma \ref{lem:Lass} gives
\begin{eqnarray}\label{eq:sum-renversee}
\sum_{G\in \mC_n} w(G)~=~  \sum_{G\in \mC_n} (-1)^{e(G)} \Vol(\Pi_G) 
~=~\frac{1}{n!} \sum_{\hh\in \ZZ^n,\si \in \Si_n,G\in\mC_n~/~ \pi(\hh,\si)\subseteq \Pi_G} (-1)^{e(G)}.
\end{eqnarray}

Let $\si$ be a permutation of $[n]$. For any vector $\hh=(h_1,\ldots,h_n)$ in $\ZZ^n$, we denote $\si(\hh)=(h_{\si(1)},\ldots,h_{\si(n)})$. For any graph $G$ labelled on $\{0,\ldots,n\}$, we denote by  $\si(G)$ the graph where the label $i$ is replaced by $\sigma(i)$ for all $i=1,\ldots,n$. 
\begin{lemma}\label{lem:desetiquetage}
Let $\hh$ be a vector in $\ZZ^n$, let $\si$ be a permutation of $[n]$ and let $G$ be a graph. Then,  $\pi(\hh,\si)\subseteq \Pi_G$ if and only if $\pi(\si^{-1}(\hh),\id)\subseteq \Pi_{\si(G)}$, where $\id$ is the identity permutation.
\end{lemma}


We omit the proof of Lemma \ref{lem:desetiquetage} which is straightforward. From this Lemma, one gets for any permutation $\sigma$ of $[n]$,
$$
\sum_{\hh\in \ZZ^n,~ G \in \mC_n \atop \pi(\hh,\si)\subseteq \Pi_G} \!\!\!\!\!(-1)^{e(G)}
= \!\sum_{\hh\in \ZZ^n,~ G\in \mC_n \atop \pi(\si^{-1}(\hh),\id)\subseteq \Pi_{\si(G)} }  \!\!\!\!\!\!\!\!(-1)^{e(G)}
= \!\sum_{\hh\in \ZZ^n,~ G\in \mC_n \atop \pi(\hh,\id)\subseteq \Pi_{G}}  \!\!\!\!\!(-1)^{e(\si^{-1}(G))}
= \!\sum_{\hh\in \ZZ^n,~ G\in \mC_n \atop \pi(\hh,\id)\subseteq \Pi_{G}}  \!\!\!\!\!(-1)^{e(G)}.
$$
where the second equality is obtained by changing the order of summations on the graphs $G$ and on the vectors $\hh$. 
Therefore, continuing Equation~\Ref{eq:sum-renversee} gives
\begin{eqnarray}\label{eq:sum-desetiquetee}
\sum_{G\in \mC_n} w(G)~=~\frac{1}{n!}\sum_{\si \in \Si_n} \sum_{\hh \in \ZZ^n,G \in \mC_n \atop \pi(\hh,\si)\subseteq \Pi_G} (-1)^{e(G)}~=~ \sum_{\hh \in \ZZ^n, G \in \mC_n \atop \pi(\hh,\id)\subseteq \Pi_G} (-1)^{e(G)}.
\end{eqnarray}

\bigskip

\titre{Killing involution.} Let  $\hh$ be a vector in $\ZZ^n$. We will now evaluate the sum $\displaystyle \!\! \sum_{G\in\mC_n~/~ \pi(\hh,\id)\subseteq \Pi_G}\!\!\!\!\!\!\!\!\!\!(-1)^{e(G)}$ thanks to a \emph{killing involution} similar to the one defined in Section~\ref{section:discrete-gas}. To the vector $\hh=(h_1,\ldots,h_n)$ we associate the \emph{centroid} $\hhb=(\hb_0,\hb_1,\ldots,\hb_n)$, where $\hb_i=h_i+\frac{i}{n+1}$ for $i=0,\ldots,n$ with the convention that $h_0=0$.  We also denotes by $G_{\hh}$ the graph on $\{0,\ldots,n\}$ whose edges are the pairs $(i,j)$ such that $|\hb_i-\hb_j|<1$.  Observe that for any graph $G$,  $\pi(\hh,\id)\subseteq \Pi_G$ if and only if $(\hb_1,\ldots,\hb_n)\in\Pi_G$ if and only if $G\subseteq G_{\hh}$.\\

We order the edges $e=(i,j)$ of the graph $G_{\hh}$ by the lexicographic order on the corresponding pairs $(\hhb_i,\hhb_j)$, that is,  $(i,j)<(k,l)$ if either $\min(\hhb_i,\hhb_j)<\min(\hhb_k,\hhb_l)$ or $\min(\hhb_i,\hhb_j)=\min(\hhb_k,\hhb_l)$ and $\max(\hhb_i,\hhb_j)<\max(\hhb_k,\hhb_l)$. For a graph $G\subseteq G_{\hh}$ and an edge $e=(i,j)$ in $G_\hh$, we denote by  $G^{>e}$ the set of edges in $G$ that are greater than $e$ and we say that  $e$ is \emph{$(G,\hh)$-active} if there is a path in $G^{>e}$ connecting $i$ and $j$. We also denote by $e^*_{G,\hh}$ the least $(G,\hh)$-active edge if there is any. We then define a mapping $\Psi_\hh$ on the set of connected spanning subgraphs of $G_{\hh}$ by: $\Psi_\hh(G)=G$ if there is no $G$-active edges and $\Psi_\hh(G)=G\oplus e^*_{G,\hh}$  otherwise.

\begin{lemma}\label{lem:killing-involution-continuum}
For any vector $\hh$ in $\ZZ^n$, the mapping $\Psi_\hh$ is an involution on the connected subgraphs of $G_\hh$.
\end{lemma}
 
The proof of Lemma~\ref{lem:killing-involution-continuum} is identical to the proof of Lemma~\ref{lem:killing-involution}. As a consequence, one gets 
$$\sum_{G\subseteq G_\hh \connect} (-1)^{e(G)}=\sum_{G\subseteq G_\hh \connect,~ \Psi_\hh(G)=G} (-1)^{e(G)}.$$

We now characterise the fixed points of $\Psi_\hh$. Let $i_0$ be the index of the least coordinate of the centroid $\hhb$ (that is, $\hb_{i_0}=\min_{i\in\{0,\ldots,n\}}(\hb_i)$).  A tree on $\{0,\ldots,n\}$ is said \emph{$\hh$-increasing} if the labels $i_0,i_1,\ldots,i_k$ of the vertices along any simple path starting at vertex $i_0$ are such that  $\hb_{i_0}<\hb_{i_1}<\cdots<\hb_{i_k}$.
\begin{lemma}\label{lem:survivors-continuum}
Let $\hh$ be a vector in $\ZZ^n$. A connected graph $G\subseteq G_\hh$ is an $\hh$-increasing tree if and only if there is no $(G,\hh)$-active edge.
\end{lemma}

\dem 
\ite We suppose first that $G$ is an $\hh$-increasing tree.  Since $G$ has no cycle, no edge in $G$ is active.  Consider now an edge $e=(i,j)\notin G$ and the nearest common ancestor $k$ of $i$ and $j$ (the root vertex of $G$ being the vertex $i_0$).  Let also $e'=(k,l)$  be an edge containing $k$ on the path in $G$ between $i$ and $j$. Since $G$ is an $\hh$-increasing tree, $\hhb_k\leq \min(\hhb_i,\hhb_j)$ and $\hhb_l\leq \max(\hhb_i,\hhb_j)$. Thus, $e'=(k,l)<e=(i,j)$ and $e$ is not $(G,\hh)$-active.\\
\ite Suppose now that there is no $(G,\hh)$-active edge.  First observe that $G$ is a tree since if $G$ had a cycle then the minimal edge in this cycle would be active. We now want to prove that the tree $G$ is $\hh$-increasing. Suppose the contrary and consider a sequence of labels $i_0,i_1,\ldots,i_r,i_{r+1}$ such that $\hb_{i_0}<\cdots<\hb_{i_{r-1}}<\hb_{i_r},\hb_{i_{r+1}}$ on a path of $G$ starting from the vertex $i_0$. Then, the edge $(i_{r-1},i_{r+1})$ belongs to $G_\hh$ (since $|\hb_{i_{r-1}}-\hb_{i_{r+1}}|<\max(|\hb_{i_{r-1}}-\hb_{i_{r}}|,|\hb_{i_{r}}-\hb_{i_{r+1}}|)<1$) and is $(G,\hh)$-active. We reach a contradiction.
\findem

From Lemma~\ref{lem:survivors-continuum}, one gets for any vector $\hh\in \ZZ^n$,
$$\sum_{G\subseteq G_\hh \connect} (-1)^{e(G)}~=~(-1)^n \# ~\{\hh\textrm{-increasing trees}\}.$$
Therefore, continuing Equation~\Ref{eq:sum-desetiquetee} gives
\begin{eqnarray}\label{eq:sum-survivors}
\sum_{G\in\mC_n} w(G)=\sum_{\hh \in \ZZ^n} ~~~ \sum_{G\subseteq G_\hh   \connect} (-1)^{e(G)}=(-1)^n \sum_{\hh \in \ZZ^n}\#~\{\hh\textrm{-increasing trees}\}.
\end{eqnarray}

\bigskip

\titre{Cayley trees.} We now relate $\hh$-increasing trees and Cayley trees.
\begin{lemma}\label{lem:Cayley-tree}
Any rooted Cayley tree  on $\{0,\ldots,n\}$ with root $i_0$ is $\hh$-increasing for exactly one vector $\hh$ in $\ZZ^n$ such that $\displaystyle \hb_{i_0}=\min_{i\in\{0,\ldots,n\}}(\hb_i)$.
\end{lemma}

\dem Let $T$ be a Cayley tree rooted on  $i_0$. The tree $T$ is $\hh$-increasing with  $\hb_{i_0}=\min_{i\in\{0,\ldots,n\}}(\hb_i)$ if and only if any vertex $j\neq i_0$ satisfies $\hb_i<\hb_j<\hb_i+1$, where $i$ is the father of $j$. The condition $\hb_i<\hb_j<\hb_i+1$ holds if and only if  either $i<j$ and $h_j=h_i$ or  $j<i$ and $h_j=h_i+1$. Therefore, tree $T$ is $\hh$-increasing with  $\hb_{i_0}=\min(\hb_i)$ if and only if for all index $i=0,\ldots,n$, the difference $h_i-h_{i_0}$ is the number of \emph{descents} in the sequence of labels $i_0,i_1,\ldots,i_s=i$ along the path of $T$ from $i_0$ to $i$ (a descent is an index $r<s$ such that $i_{r+1}<i_r$). Knowing that $h_0=0$ completes the proof.
\findem

It is well known that the number of rooted Cayley trees on $\{0,\ldots,n\}$ is $(n+1)^n$. Thus, Lemma~\ref{lem:Cayley-tree} gives 
\begin{eqnarray}
\sum_{G\in \mC_n} w(G)= (-1)^n \sum_{\hh \in \ZZ^n} \# ~ \hh\textrm{-increasing trees}=(-1)^n (n+1)^n. 
\end{eqnarray}

This completes the proof of Equation~\Ref{eq:Mayer-continuum-gas-intro} and answers the question of Labelle \emph{et al.} \cite{Leroux:mayers-theory2}.\\

\bigskip

\noindent \textbf{Acknowledgements:} I would like to thank Pierre Leroux who warmly encouraged me to work on the combinatorial interpretations of identities arising from Mayer's theory of cluster integrals.

\bibliography{../../../biblio/allref}
\bibliographystyle{plain}

\end{document}